\documentclass[mnsc,nonblindrev]{informs3} % current default for manuscript submission

\OneAndAHalfSpacedXI % current default line spacing
\usepackage{hyperref}

\usepackage{endnotes}
\let\footnote=\endnote
\let\enotesize=\normalsize
\def\notesname{Endnotes}%
\def\makeenmark{$^{\theenmark}$}
\def\enoteformat{\rightskip0pt\leftskip0pt\parindent=1.75em
  \leavevmode\llap{\theenmark.\enskip}}

\usepackage{amsmath}
\usepackage{mathtools}
\usepackage{stackengine}
\usepackage{natbib}
 \bibpunct[, ]{(}{)}{,}{a}{}{,}%
 \def\bibfont{\small}%

\usepackage{algorithm}
\usepackage{algorithmic}

\usepackage{makecell}

\usepackage{ulem}
\usepackage{booktabs}
\usepackage{multirow}
\usepackage{comment}

\usepackage{undertilde_closer}

\usepackage{hyperref}

\usepackage{xcolor}
\hypersetup{
    colorlinks,
    linkcolor={red!50!black},
    citecolor={blue!50!black},
    urlcolor={blue!80!black}
}
\TheoremsNumberedThrough     % Preferred (Theorem 1, Lemma 1, Theorem 2)
\ECRepeatTheorems

\EquationsNumberedThrough    % Default: (1), (2), ...
\newenvironment{assumption'}[1]
  {\renewcommand{\theassumption}{\arabic{assumption}$'$}%
   \addtocounter{assumption}{-1}%
   \begin{assumption}}
  {\end{assumption}}

\usepackage{amsmath}
\usepackage{amsfonts}
\usepackage{amssymb}
\usepackage{dsfont}
\newenvironment{alg'}[1]
  {\renewcommand{\thealgorithm}{\ref{#1}$'$}%
   \addtocounter{algorithm}{-1}%
   \begin{breakablealgorithm}}
  {\end{breakablealgorithm}}

\newenvironment{lemma'}[1]
  {\renewcommand{\thelemma}{\ref{#1}$'$}%
   \addtocounter{lemma}{-1}%
   \begin{lemma}}
  {\end{lemma}}

\begin{document}
%%%%%%%%%%%%%%%%

%\setcounter{page}{1}
%\DOI{10.1287/xxxx.0000.0000}%

% \RUNAUTHOR{A,B and C}

\RUNTITLE{RT}

% \TITLE{Bundling Information Goods of Decreasing Value}
\TITLE{A Regeneration-based a Posteriori Error Bound for a Markov Chain Stationary Distribution Truncation Algorithm}


\ARTICLEAUTHORS{%
\AUTHOR{Peter W. Glynn}
\AFF{Stanford University}
\AUTHOR{Zeyu Zheng}
\AFF{University of California Berkeley}
} % end of the block

\ABSTRACT{%
When the state space of a discrete state space positive recurrent Markov chain is infinite or very large, it becomes necessary to truncate the state space in order to facilitate numerical computation of the stationary distribution. This paper develops a new approach for bounding the truncation error that arises when computing approximations to the stationary distribution. This rigorous a posteriori error bound exploits the regenerative structure of the chain and assumes knowledge of a Lyapunov function. Because the bound is a posteriori (and leverages the computations done to calculate the stationary distribution itself), it tends to be much tighter than a priori bounds. The bound decomposes the regenerative cycle into a random number of excursions from a set $K$ defined in terms of the Lyapunov function into the complement of the truncation set $A$. The bound can be easily computed, and does not (for example) involve a linear program, as do some other error bounds.
}%

%\KEYWORDS{deterministic inventory theory; infinite linear programming duality;
%  existence of optimal policies; semi-Markov decision process; cyclic schedule}

% \KEYWORDS{ } 
%
% \HISTORY{This paper was first submitted on }

\maketitle
%%%%%%%%%%%%%%%%%%%%%%%%%%%%%%%%%%%%%%%%%%%%%%%%%%%%%%%%%%%%%%%%%%%%%%

% Samples of sectioning (and labeling) in OPRE
% NOTE: (1) \section and \subsection do NOT end with a period
%       (2) \subsubsection and lower need end punctuation
%       (3) capitalization is as shown (title style).
%
%\section{Introduction.}\label{intro} %%1.
%\subsection{Duality and the Classical EOQ Problem.}\label{class-EOQ} %% 1.1.
%\subsection{Outline.}\label{outline1} %% 1.2.
%\subsubsection{Cyclic Schedules for the General Deterministic SMDP.}
%  \label{cyclic-schedules} %% 1.2.1
%\section{Problem Description.}\label{problemdescription} %% 2.

% Text of your paper here

\section{Introduction}

Let $X=(X_n:n\geq 0)$ be an irreducible positive recurrent Markov chain taking values in an either countably infinite state space $S$ or a finite state space $S$ with a very large number of elements. If $P=(P(x,y):x,y\in S)$ is the one-step transition matrix, there is a unique probability distribution $\pi=(\pi(x):x\in S)$ that satisfies the linear system
\begin{align}
\pi = \pi P,\label{eq-pieqpiP}
\end{align}
where we have chosen here to encode $\pi$ as a row vector. The distribution $\pi$ is called the \textit{stationary distribution} of $X$ (or, equivalently, the equilibrium distribution, steady-state distribution, or invariant distribution of $X$). The stationary distribution describes the long-run behavior of $X$, and is therefore of interest in many applications. If we wish to solve \eqref{eq-pieqpiP} numerically, this linear system must clearly be truncated to a finite linear system.

This truncation problem dates back to \citet{senetaFiniteApproximationsInfinite1967}. Much of the associated literature on the truncation problem focuses only on establishing convergence of various proposed truncation algorithms, and no error bounds are provided. However,   \citet{tweedieTruncationApproximationsInvariant1998} developed a priori error bounds, under the assumption of geometrically ergodicity, which can be calculated prior to numerically computing the truncated approximation. Further such a priori error bounds have been developed by exploiting different structural assumptions on the Markov chain; see, for example, \citet{herveApproximatingMarkovChains2014}, \citet{masuyamaErrorBoundsAugmented2015}, and \citet{liuErrorBoundsAugmented2018}. 
% Zeyu, are these a priori error bounds? 
More recently, \citet{kuntzBoundingStationaryDistributions2019} developed an approach to computing a posteriori error bounds that involves solving a linear program with the number of decision variables equal to the number of states in the truncation set. \citet{kuntzStationaryDistributionsContinuoustime2021} surveys the truncation literature, providing an overview of both truncation algorithms and more recent efforts to develop computable error bounds.
% Note Peter first put Ledoux et al. 2018 but I think he meant Liu et al. check it.

To discuss our approach, we let $P_x(\cdot)=P(\cdot|X_0=x)$ and $E_x(\cdot)=E[\cdot|X_0=x]$ be the probability and expectation on the path-space of $X$ associated with conditioning on $X_0=x\in S$. Our approach relies on the well known fact that for any $z\in S$ (known henceforth as the \textit{regeneration state}), we may write $\pi$ as
\begin{align}
\pi(\cdot) &=\frac{E_z\sum_{j=0}^{\tau(z)-1}I(X_j\in \cdot)}{E_z \tau(z)},\label{eq_12_regenerative_representation}
\end{align}
where $\tau(z)=\inf\{n\geq 1:X_n=z\}$ is the first return time to $z$ and $I(A)$ is the indicator random variable associated with the event $A$. The equality \eqref{eq_12_regenerative_representation} asserts that the stationary distribution may be expressed in terms of the ratio of two expectations, each involving the behavior over a regenerative ``$z$-cycle''; see \citet{chungMarkovChains1967} for details. In view of \eqref{eq_12_regenerative_representation}, our approach involves approximating the numerator and denominator of \eqref{eq_12_regenerative_representation} separately.

To obtain error bounds for our truncation, we construct upper and lower bounds on the two expectations arising in \eqref{eq_12_regenerative_representation}. It turns out that lower bounds can be easily obtained. However, our upper bounds require knowledge of a Lyapunov bound that holds on the complement of some finite subset $K \supseteq \{z\}$. In this paper, we develop a new approach to constructing an upper bound to the expectations appearing in \eqref{eq_12_regenerative_representation}. This approach is different from the error bounds developed in Infanger and Glynn (2022). The error bounds there exploit a ratio formula, involving cycles defined by returns to  $K$, rather than the singleton $\{z\}$ used in this paper. Furthermore, those error bounds are formulated in terms of a linear program or the solution to a finite state Poisson's equation, neither of which arises in computing the error bound found in the current paper. 

The bound developed in this paper splits the regenerative $z$-cycles into a random number of excursions from $z$ into the complement of the truncation set and back into $K$; see Theorem 1. We can then use this result to obtain a bound on the (weighted) total variation distance between $\pi$ and its truncated approximation $\tilde\pi$; see Theorem 2. 

Section 2 of this paper is devoted to developing our bounds and proving our key theoretical results. Section 3 reports on our computational experience with these new bounds.

\section{The Regenerative Bounds}

Fix $z\in S$ and put $S'=S-\{z\}$. Let $\nu=(\nu(x):x\in S')$, $e=(e(x):x\in S')$, and $B=(B(x,y):x,y\in S')$ have entries given by $\nu(x)=P(z,x)$, $e(x)=1$ and $B(x,y)=P(x,y)$ for $x,y\in S'$. If we encode $\nu$ as a row vector and $e$ as a column vector, then \eqref{eq_12_regenerative_representation} asserts that
\begin{align}
\pi(x)=\frac{(\nu\sum_{j=0}^{\infty}B^j)(x)}{1+\nu\sum_{j=0}^{\infty}B^j e}\label{eq21}
\end{align}
for $x\in S'$, with 
\begin{align}
\pi(z)=(1+\nu\sum_{j=0}^{\infty}B^je)^{-1}.\label{eq22}
\end{align}
Let $A\subseteq S$ be a finite subset containing $z$ that represents the largest possible truncation set that is computationally feasible within the computing budget that is available. Given $A$, we put $A'=A-\{z\}$, and let $\tilde\nu=(\tilde \nu(x):x\in A'), \tilde e=(\tilde e(x):x\in A')$, $\tilde B=(\tilde B(x,y):x,y\in  A')$ and $\tilde{p}=(\tilde{p}(x):x\in A')$ have entries given by $\tilde \nu(x)=P(z,x)$, $\tilde e(x)=1$,  $\tilde B(x,y)=P(x,y)$ and $\tilde{p}(x) = P(x,z)$ for $x,y\in A'$. Since $P$ is irreducible, there exists a finite length path leading to a state in $\{z\}\cup A^c$ having positive probability. It follows that $\tilde B^n(x,y)\rightarrow 0$ for each $x,y\in {A'}$, so that $I-\tilde B$ is non-singular; see Lemma B.1. in \citet{senetaNonnegativeMatricesMarkov2006}. If we now encode $\tilde \nu$ as a row vector and $\tilde e$ as a column vector, our truncation approximation to $\pi$ is then given by 
\begin{align}
\tilde \pi(x) = \frac{(\tilde \nu(I-\tilde B)^{-1})(x)}{1+\tilde \nu(I-\tilde B)^{-1}\tilde e}\label{eq23}
\end{align}
for $x\in {A'}$, with 
\begin{align}
\tilde \pi(y)= \begin{cases} (1+\tilde \nu(I-\tilde B)^{-1}\tilde e)^{-1}, & y=z,\\
0, & y\in A^c.
\end{cases}\label{eq24}
\end{align}

Consider a sequence of truncation sets $(A_n:n\geq 1)$ such that $z\in A_1$ with $A_n\nearrow S$ as $n\rightarrow\infty$, and suppose that $(\tilde \pi_n:n\geq 1)$ is the sequence of corresponding truncation approximations.

Our first result asserts that for each $x\in S$, $\tilde \pi_n(x)\rightarrow \pi(x)$ as $n\rightarrow\infty$. 

\begin{proposition} Suppose that $X$ is irreducible and positive recurrent. Then, for each $x\in S$, $\tilde \pi_n(x)\rightarrow \pi(x)$ as $n\rightarrow\infty$. 
\end{proposition}
\textit{Proof.} Note that $(I-\tilde B)^{-1}=\sum_{j=0}^{\infty}\tilde B^j$, so for $y\in {A'}$, 
\begin{align*}
(\tilde \nu &\sum_{j=0}^{\infty}\tilde B^j)(y) = \sum_{x\in A'}^{}P(z,x)\sum_{j=0}^{\infty}\tilde B^j(x,y)\\
&=\sum_{x\in {A'}}^{}P(z,x)(1+\tilde B(x,y)+\sum_{j=2}^{\infty} \sum_{\substack{w_1,...,w_{j-1}\\w_i\in A'}}^{}\tilde B(x,w_1)\tilde B(w_1,w_2)... \tilde B(w_{j-1},y)\\
&= \sum_{x\in {A'}}^{}P(z,x)(1+P(x,y)+\sum_{j=2}^{\infty}\sum_{\substack{w_1,...,w_{j-1}\\w_i\in A'}}^{}P(x,w_1)P(w_1,w_2)...P(w_{j-1}, y))\\
&=\sum_{j=1}^{\infty} P_z(X_j=y, X_1\not\in A^c,..., X_j\not\in A^c, \tau(z)>j)\\
&=\sum_{j=0}^{\infty}P_z(X_j=y, T\land \tau(z)>j),
\end{align*}
where $T=\inf\{n\geq 0: X_n\in A^c\}$ is the first hitting time of $A^c$. Hence, we conclude that
\begin{align*}
(\tilde\nu \sum_{j=0}^{\infty}\tilde B^j\tilde e_z)(y) = E \smashoperator{\sum_{j=0}^{(\tau(z)\land T)-1}}I(X_j=y)
\end{align*}
for $y\in A'$. It follows that for $n$ large enough that $y\in A_n$ (with $y\neq z$),
\begin{align}
\tilde \pi_n(y) = E_z\smashoperator{\sum_{j=0}^{(\tau(z)\land T_n)-1}}I(X_j=y)/E_z(\tau(z)\land T_n)\label{eq25}
\end{align}
where $T_n$ is the first hitting time of $A_n$. The Monotone Convergence Theorem, applied separately to the numerator and denominator of \eqref{eq25}, then guarantees that $\tilde \pi_n(y)\rightarrow \pi(y)$ as $n\rightarrow\infty$ for $y\neq z$. A similar argument ensures that $\tilde\pi_n(z)\rightarrow\pi(z)$ as $n\rightarrow\infty$.
\Halmos

\begin{remark}
The probability $\tilde{\pi}_n$ is identical to the stationary distribution of the Markov chain that transitions according to  $P$ when inside $A_n$, and transitions immediately to $z$ whenever the chain attempts to transition to $A_n$'s. Consequently, $\tilde{\pi}_n$ corresponds to what is known in the literature as ``fixed state" transition-augmentation. Proposition 1 has been previously established for fixed state schemes by \citet{wolf1980approximation}.
\end{remark}

\begin{remark} Note that if $X$ is irreducible and null recurrent, the unique (up to a multiplicative constant) non-negative solution to \eqref{eq-pieqpiP} has the property that for $x,y\in S$,
\begin{align}
\frac{\pi(x)}{\pi(y)} = \frac{E_z\sum_{j=0}^{\tau(z)-1}I(X_j=x)}{E_z\sum_{j=0}^{\tau(z)-1} I(X_j=y)},\label{eq26}
\end{align}
where the numerator and denominator of the right-hand side of \eqref{eq26} are both guaranteed to be positive and finite; see \citet{chungMarkovChains1967}. It follows that when $(A_n:n\geq 1)$ satisfies A1, then
\begin{align*}
\frac{\tilde\pi_n(x)}{\tilde\pi_n(y)}\rightarrow \frac{\pi(x)}{\pi(y)}
\end{align*}
as $n\rightarrow\infty$, so that our truncation algorithm is convergent even when $X$ is null recurrent.
\end{remark}

We now proceed to developing error bounds for our approximation. Given a non-negative function $r: S\rightarrow \mathbb{R}_+$, our goal is to develop upper and lower bounds of $\alpha(r)\triangleq \sum_x \pi(x)r(x)$. We accomplish this goal by recognizing that \eqref{eq_12_regenerative_representation} implies that
\begin{align}
\alpha(r)= \frac{r(z)+\nu\sum_{j=0}^{\infty}B^jr}{1+\nu\sum_{j=0}^{\infty}B^je},\label{eq27}
\end{align}
where $r=(r(x):x\in S')$ is the column vector in which the $x$'th (for $x\in S'$) entry is given by $r(x)$. Suppose that we can find lower bounds $\utilde \kappa (r)$ and $\utilde \kappa(e)$ and upper bounds $\tilde \kappa(r)$ and $\tilde \kappa(e)$ for which
\begin{align*}
\utilde\kappa(r) \leq r(z) + \sum_{j=0}^{\infty}\nu B^j r \leq \tilde \kappa(r),
\end{align*}
and 
\begin{align*}
\utilde \kappa(e)\leq 1+\sum_{j=0}^{\infty}\nu B^je\leq \tilde\kappa(e).
\end{align*}
With these bounds, we note that
\begin{align}
\frac{\utilde \kappa(r)}{\tilde \kappa(e)}\leq \alpha(r)\leq \frac{\tilde \kappa(r)}{\utilde\kappa(e)},\label{eq28}
\end{align}
thereby providing our desired upper and lower bounds on $\alpha(r)$.

A lower bound on $\utilde \kappa(r)$ is easy to develop. Let $\tilde{r}=(\tilde{r}(x):x\in{A'})$ be such that $\tilde{r}(x) = r(x)$ for $x\in{A'}$. Then, the proof of Proposition 1 shows that
\begin{align*}
r(z) + \tilde \nu\sum_{j=0}^{\infty}\tilde B^j\tilde{r} &= E_z\sum_{j=0}^{(\tau(z)\land T)-1}r(X_j)\\
&\leq E_z \sum_{j=0}^{\tau(z)-1}r(X_j)\\
&= r(x)+\nu \sum_{j=0}^{\infty}B^jr,
\end{align*}
so that we may use
\begin{align*}
r(z) + \tilde \nu (I-\tilde B)^{-1}\tilde r \overset{\Delta}{=}\utilde \kappa(r)
\end{align*}
as our lower bound. Similarly
\begin{align*}
\utilde \kappa(e) \overset{\Delta}{=} 1+ \tilde \nu(I-\tilde B)^{-1}\tilde e\leq 1+\nu_z\sum_{j=0}^{\infty}B^j e.
\end{align*}
Note that $\tilde \pi r$ can be represented as $\utilde\kappa(r)/\utilde\kappa(e)$. It follows that
\begin{align}
\frac{\utilde \kappa(r)}{\tilde \kappa(e)}\leq \tilde \pi r \leq \frac{\tilde \kappa (r)}{\utilde \kappa(e)}\label{eq29}
\end{align}
and hence, in view of \eqref{eq28}, we have
\begin{align}
|\pi r-\tilde \pi r| \leq \frac{\tilde \kappa(r)}{\utilde \kappa(e)} - \frac{\utilde \kappa(r)}{\tilde \kappa(e)}. \label{eq210}
\end{align}
We now turn to developing an upper bound $\tilde \kappa(r)$. We will assume the existence of \textit{Lyapunov functions} $g_1$ and $g_2$ satisfying the conditions below.\\

\begin{assumption}\label{A1} Suppose there exists a finite set $K\subseteq A$ and containing $z$, for which there exists known non-negative functions $g_1,g_2, h_1,h_2$ such that
\begin{align}
\sum_{y\in K^c}^{}P(x,y)g_1(y) \leq g_1(x)-r(x),\label{eq211}\\
\sum_{y\in K^c}^{}P(x,y)g_2(y) \leq g_2(x)-1,\label{eq212}
\end{align}
for $x\in K^c$, with
\begin{align*}
\sum_{y\in A^c}P(x,y) g_1(y)\leq h_1(x),\\
\sum_{y\in A^c}^{}P(x,y)g_2(y) \leq h_2(x)
\end{align*}
for $x\in A$.\\
\end{assumption}
It is well known that \eqref{eq211} and \eqref{eq212} respectively imply
\begin{align}
E_x\sum_{j=0}^{T-1}r(X_j) \leq g_1(x)\label{eq213}
\end{align}
and 
\begin{align}
E_xT \leq g_2(x)\label{eq214}
\end{align}
for $x\in K^c$; see p.344 of \citet{meynMarkovChainsStochastic2012}. With these bounds on excursions of $X$ from $A^c$ into $K$ in place, we can now use Assumption \ref{A1} to develop a computable upper bound $\tilde\kappa(r)$. 

To accomplish this, we put $\Gamma_0=0$, and $T_1=T$. For $i\geq 1$, put
\begin{align*}
\Gamma_i = \inf\{n>T_i:X_n\in K\}
\end{align*}
and
\begin{align*}
T_{i+1} = \inf\{n>\Gamma_i: X_n\in A^c\}.
\end{align*}
The recurrence of $X$ implies that $\Gamma_i<\infty$ a.s. and $T_i<\infty$ a.s. for $i\geq 1$. Put
\begin{align*}
N=\inf\{n\geq 1: \tau(z)<\Gamma_n\}
\end{align*}
and observe that
\begin{align}
\sum_{j=0}^{\tau(z)-1}r(X_j) &= \sum_{i=1}^{N}\sum_{j=\Gamma_{i-1}}^{(\Gamma_i\land \tau(z))-1}r(X_j)\nonumber\\
&= \sum_{i=1}^{\infty} I(N\geq i) \sum_{j=\Gamma_{i-1}}^{(\Gamma_i\land \tau(z))-1}r(X_j)\label{eq215}.
\end{align}
Note that for $x\in {A'}$, \eqref{eq213} implies that
\begin{align}
E_x&\smashoperator{\sum_{j=0}^{(\Gamma_1\land \tau(z))-1}} r(X_j) = E_x\smashoperator{\sum_{j=0}^{(T_1\land \tau(z))-1}}r(X_j) + E_x I(\tau(z)>T_1)\sum_{j=T_1}^{\Gamma_1-1}r(X_j)\nonumber\\
&= \sum_{j=0}^{\infty} E_x r(X_j)I(\tau(z)\land T_1>j) + \sum_{j=0}^{\infty}E_xI(\tau(z)\land T_1 >j, T_1 = {j+1})\sum_{k=j+1}^{\Gamma_1-1}r(X_k)\nonumber\\
&\leq \sum_{j=0}^{\infty}(\tilde B^j\tilde r )(x) + \sum_{j=0}^{\infty} E_xI(\tau(z)\land T_1>j)I(T_1 = j+1) g_1(X_{j+1})\nonumber\\
&=((I-\tilde B)^{-1} \tilde r )(x)+\sum_{j=0}^{\infty}\sum_{w\in A'}^{}P_x(\tau(z) \land T_1>j, X_j=w)\sum_{y\in A^c}^{}P(w,y) g_1(y)\nonumber\\
&\leq ((I-\tilde B)^{-1} \tilde r )(x)+\sum_{j=0}^{\infty}\sum_{w\in A'}^{}P_x(\tau(z)\land T_1>j, X_j=w)h_1(w)\nonumber\\
&=((I-\tilde B)^{-1} \tilde r )(x) + \sum_{j=0}^{\infty}E_xh_1(X_j)I(\tau(z)\land T_1>j)\nonumber\\
&=((I-\tilde B)^{-1}\tilde r )(x) + ((I-\tilde B)^{-1} \tilde h_{1})(x)\label{eq216},
\end{align}
where $\tilde h_{i}=(\tilde h_{i}(x):x\in {A'})$ is the column vector for which  $\tilde h_{i}(x)=h_i(x)$ for $x\in {A'},i=1,2$. 

Set ${K'}=K-\{z\}$. For $i\geq 2$ and $x\in K$, \eqref{eq216} implies that
\begin{align}
E_xI(N\geq i) \smashoperator{\sum_{j=\Gamma_{i-1}}^{(\Gamma_i\land\tau(z))-1}}r(X_j) & = E_x I(\tau(z)>\Gamma_{i-1})E_x\Big[\sum_{j=\Gamma_{i-1}}^{(\Gamma_i\land \tau(z))-1} r(X_j)|X_0,X_1,...,X_{\Gamma_{i-1}}\Big]\nonumber\\
&\leq E_x I(\tau(z)>\Gamma_{i-1}) ((I-\tilde B)^{-1})(\tilde r+\tilde h_{1})(X_{\Gamma_{i-1}})\nonumber\\
&\leq \|{(I-\tilde B)^{-1}(\tilde r+\tilde h_{1})}\|_{{K'}}\cdot P_x(\tau(z)>\Gamma_{i-1})\label{eq217}
\end{align}
where
\begin{align*}
\|{v}\|\overset{\Delta}{=}\max\{|v(x)|:x\in {K'}\}
\end{align*}
for a generic column vector $v=(v(x):x\in {K'})$. 

On the other hand, for $x\in {K'}$ and $i\geq 1$,
\begin{align}
P_x(\tau(z)>\Gamma_i) & = P_x(\tau(z)>\Gamma_{i-1}) - P_x(\Gamma_{i-1}<\tau(z)\leq \Gamma_{i})\nonumber\\
&\leq P_x(\tau(z)>\Gamma_{i-1}) - P_x(\Gamma_{i-1}<\tau(z)<T_i)\nonumber\\
&= P_x\left(\tau(z)>\Gamma_{i-1}\right) - \sum_{j=0}^{\infty}P_x\left(\Gamma_{i-1}<\tau(z), \tau(z)\land T_i>\Gamma_{i-1}+j, X_{\Gamma_{i-1}+j+1}=z\right)\nonumber\\
&=P_x\left(\tau(z)>\Gamma_{i-1}\right) - \sum_{j=0}^{\infty} E_x I(\Gamma_{i-1}<\tau(z)) \sum_{y\in {A'}}^{}\tilde B^j(X_{\Gamma_{i-1}},y)P(y,z)\nonumber\\
&=P_x(\tau(z)>\Gamma_{i-1}) - \sum_{y\in {A'}}^{}E_xI(\tau(z)>\Gamma_{i-1})((I-\tilde B)^{-1}(X_{\Gamma_{i-1}}, y))P(y,z)\nonumber\\
&\leq P_x(\tau(z)>\Gamma_{i-1})(1-\delta),\label{eq218}
\end{align}
where
\begin{align*}
\delta \overset{\Delta}{=} \min_{x\in {K'}} \sum_{y\in {A'}}^{}(I-\tilde B)^{-1}(x,y)P(y,z) = \min_{x\in {K'}}((I-\tilde{B})^{-1}\tilde{p})(x).
\end{align*}
It follows from \eqref{eq217} and \eqref{eq218} that for $i\geq 1$,
\begin{align}
P_z(\tau(z)>\Gamma_i) \leq (1-\beta)(1-\delta)^{i-1},\label{eq219}
\end{align}
where $\beta = P_z(\tau(z)<T_1) = P(z,z)+\tilde{\nu}(I-\tilde B)^{-1}\tilde{p}$. Combining \eqref{eq215}, \eqref{eq216}, \eqref{eq219} (and noting that a similar argument applies to $e$ in place of $r$), we arrive at the following theorem.

\begin{theorem} \label{thm1} Assume Assumption \ref{A1} and suppose that $\delta>0$. Then,
\begin{align*}
E_z \sum_{j=0}^{\tau(z)-1}r(X_j) \leq \tilde \kappa(r)
\end{align*}
and 
\begin{align*}
E_z \tau(z) \leq \tilde\kappa(e),
\end{align*}
where 
\begin{align*}
\tilde \kappa(r) = r(z) + \tilde \nu(I-\tilde B)^{-1}(\tilde r + \tilde h _{1}) + \frac{1-\beta}{\delta}\|{(I-\tilde B)^{-1}(\tilde r + \tilde h_{1})}\|_{{K'}}
\end{align*}
and
\begin{align*}
\tilde\kappa(e) = 1+\tilde\nu(I-\tilde B)^{-1}(\tilde e+\tilde h_{2})+ \frac{1-\beta}{\delta}\|{(I-\tilde B)^{-1}(\tilde e+\tilde h_{2})}\|_{{K'}}.
\end{align*}
\end{theorem}

With Theorem \ref{thm1} in hand, we can now compute the upper bound of \eqref{eq210} on $|\pi r-\tilde\pi r|$. This leads to the inequality 
\begin{align}
|\pi r-\tilde \pi r|&\leq \frac{\tilde\kappa(r)\tilde\kappa(e)-\utilde\kappa(r)\utilde\kappa(e)}{\utilde\kappa(e)\tilde\kappa(e)}\nonumber\\
& = \frac{\utilde\kappa(r)\Delta_2 + \utilde\kappa(e)\Delta_1+\Delta_1\Delta_2}{\utilde\kappa(e)\tilde\kappa(e)},\label{eq220}
\end{align}
where
\begin{align*}
\Delta_1=\tilde\nu(I-\tilde B)^{-1}\tilde h_{1} + \frac{1-\beta}{\delta}\|{(I-\tilde B)^{-1}(\tilde r+\tilde h_{1})}\|_{{K'}},\\
\Delta_2= \tilde \nu(I-\tilde B)^{-1} \tilde h_{2} + \frac{1-\beta}{\delta}\|{(I-\tilde B)^{-1}(\tilde e + \tilde h_{2})}\|_{{K'}}.
\end{align*}

We can also use \eqref{eq220} to obtain a bound on the total variation distance between $\pi$ and $\tilde \pi$. In particular, we consider the $r$-weighted total variation distance given by
\begin{align*}
\|{\pi-\tilde \pi}\|_{r}= \sup\{|\pi w-\tilde\pi w|:|w|\leq r\}.
\end{align*}
For a given function $w:S\rightarrow\mathbb R$, we can write $w=w_+-w_-$, where $w_+(x)=\max(w(x),0)$ and $w_-(x)=\max(-w(x),0)$. Then,
\begin{align*}
|(\pi-\tilde \pi)w|\leq |(\pi-\tilde\pi)w_+| + |(\pi -\tilde \pi)w_-|.
\end{align*}
So
\begin{align*}
\| \pi-\tilde\pi \|_{r}\leq 2\sup\{|\pi w-\tilde \pi w|:0\leq w\leq r\}.
\end{align*}
For $w$ such that $0\leq w\leq r$, condition \ref{A1} holds also with $w$ replacing $r$, and with exactly the same choice of $g_1$ and $h_1$. If we set $w_z=(w_z(x):x\in S')$ and $\tilde w_z(x)=(w_z(x):x\in A')$ with $w_z(x)=w(x)$ for $x\in S'$ and $\tilde w_z(x)=w(x)$ for $x\in A'$, Theorem  \ref{thm1} establishes that 
\begin{align*}
\utilde \kappa(w) \leq w(z) + \sum_{j=0}^{\infty}B^jw_z \leq \tilde \kappa(w),
\end{align*}
where 
\begin{align*}
\utilde \kappa(w) &= w(z) + \tilde \nu (I-\tilde B)^{-1}\tilde w,\\
\tilde \kappa(w) &= \utilde \kappa(w) + \tilde \nu (I-\tilde B)^{-1} \tilde h_{1} + \frac{1-\beta}{\delta}\|{(I-\tilde B)^{-1}(\tilde w_z + \tilde h_{1})}\|_{{K'}}.
\end{align*}
Because $0\leq w \leq r$ and $(I-\tilde B)^{-1}$ is non-negative, it is evident that 
\begin{align*}
\|{(I-\tilde B)^{-1}(\tilde w_z + \tilde h_{1})}\|_{{K'}} \leq \|{(I-\tilde B)^{-1}(\tilde r+\tilde h_{1})}\|_{{K'}}.
\end{align*}
Consequently, \eqref{eq220} translates (for $0\leq w \leq r$) into the bound
\begin{align*}
|\pi w-\tilde \pi w| &\leq \frac{\utilde \kappa(w)\Delta_2 + \utilde \kappa(e)\Delta_1 + \Delta_1\Delta_2}{\utilde\kappa(e)\tilde\kappa(e)}\\
&\leq \frac{\utilde\kappa(r)\Delta_2 + \utilde\kappa(e)\Delta_1 + \Delta_1\Delta_2}{\utilde\kappa(e)\tilde\kappa(e)}.
\end{align*}
We have therefore obtained the following total variation bound on the distance between $\pi$ and $\tilde \pi$.

\begin{theorem}\label{thm2} Suppose Assumption \ref{A1} holds and suppose that $\delta>0$. Then,
\begin{align*}
\|{\pi-\tilde\pi}\|_{r}\leq 2 \frac{(\utilde\kappa(r)\Delta_2 + \utilde\kappa(e)\Delta_1 + \Delta_1\Delta_2}{\utilde\kappa(e)\tilde\kappa(e)}.
\end{align*}
\end{theorem}

\begin{remark} In most applications, $r$ can be chosen greater than or equal to $e$. In that case, $g_1$ and $h_1$ can then also play the role of $g_2$ and $h_2$ in Assumption \ref{A1}. We state Assumption \ref{A1} with two different choices (i.e. $(g_1,h_1)$ and $(g_2,h_2)$), so that one can customize the Lyapunov function for $e$, thereby potentially obtaining a tighter upper bound $\tilde\kappa(e)$.
\end{remark}

\section{Numerical experiment}
In this section, we consider the computation of the error bounds for stationary expectations developed in this paper. In particular, we compare our proposed regenerative algorithm with the truncation algorithm in \cite{infanger2022new}.  We consider two models. 
% Note that the two algorithms requires to compute $\tilde{\kappa}$ and $\utilde{\kappa}$. Hence, the total computating time is the sum of the time of computing these 2 quantities.
% \begin{algorithm}[h]
%         \caption{Regenerative Bounds Computation}
%         \begin{algorithmic}[1]
%         \STATE Given $S,P,z\in S, A, K, r, g_1,g_2,h_1,h_2$
%         \STATE Define $\tilde{A} = A-\left\{z\right\},\tilde{B}(x,y) = \left. P(x,y)\right\vert_{\tilde{A}}  $
        
%         \STATE Compute 

%         $$\beta = P(z,z) + \tilde{\nu} (I-\tilde{B})^{-1} \tilde{p}$$
%         $$\delta = \min_{x\in {K'}} \left((I-\tilde{B})^{-1} \tilde{p}\right)(x) $$
%         $$\tilde{\kappa}(r) = r(z) + \tilde{\nu} (I-\tilde{B})^{-1} \left(\tilde{r} + \tilde{h_1}\right) +\frac{1-\beta}{\delta}\left\lVert (I-\tilde{B})^{-1} (\tilde{r} + \tilde{h_1})\right\rVert _{{K'}}  $$
%         $$\tilde{\kappa}(e) = 1 + \tilde{\nu} (I-\tilde{B})^{-1} \left(\tilde{e} + \tilde{h_2}\right) +\frac{1-\beta}{\delta}\left\lVert (I-\tilde{B})^{-1} (\tilde{e} + \tilde{h_2})\right\rVert _{{K'}}$$
%         $$\utilde{\kappa} (r) = r(z) + \tilde{\nu} (I-\tilde{B})^{-1} \tilde{r}$$
%         $$\utilde{\kappa} (e) = 1 + \tilde{\nu} (I-\tilde{B})^{-1} \tilde{e}$$

%         \STATE Return: $\tilde{\kappa}(r),\tilde{\kappa}(e),\utilde{\kappa} (r),\utilde{\kappa} (e)$
%         \end{algorithmic}
%         \end{algorithm}
\textbf{1. G/M/1 Model.}
Here, the state space $S=\{0,1,2,\ldots\}$ and the transition matrix is given by 
$$P(x,y) = \begin{cases}
    \beta_{x+y-1},~1\leq y \leq x+1\\
     0, ~y>x+1 \\
     1 - \sum_{i=0}^{x} \beta_i, ~y=0,
\end{cases}$$
where the $\beta_i$'s are given by 
$$\beta_i = \int_{0}^{2.01} \frac{e^{-t} t^i}{2.01 i!}   \,dt , i\geq 0.$$
The set $K$ was chosen as $K=\{0,1,\cdots,200\}$, $z=0$, and $A$ was selected to be of the form $A=\{0,1,\cdots,a-1\}$. We use the Lyapunov functions
$$g_1(x) = 300 x^2,~~g_2(x) = 300x,$$
$$h_i(x) = \sum_{y\in A^c} P(x,y) g_i(y) = \begin{cases}
    0,~~0\leq x \leq a-1\\
     300\beta_0  \cdot (a+1)^{3-i},~~x=a
\end{cases},\,i=1,2.$$ We compute the bounds for $\pi r$, where $r(x) = x$, as a function of $a$.

\textbf{2. Random Walk.}
We consider a random walk on $S = \left\{0,1,2,\cdots\right\}$ with $r(x) = x/2$. 
For each positive integer, $j\geq 1$, the random walk has a probability $\frac{1}{3}$ of moving to $j+1$ and a probability $\frac{2}{3}$ of moving to $j-1$.
At state $0$, the random walk has probability one of moving to state 1. We selected $K=\{0,1,\cdots,300\}$, $z=0$, $A=\{0,1,\cdots,a-1\}$, and used the Lyapunov functions  $g_1(x) = g_2(x)=x^2$,  and
% As a result, the stochastic matrix of this model is 
% $$P=\begin{bmatrix}
%     0 & 1 & & & \\
%     \frac{2}{3} & 0 & \frac{1}{3} & & \\
%     & \frac{2}{3} & 0 & \frac{1}{3} & \\
%     & & \frac{2}{3} & \cdots & \\
%     & & & & \cdots 
% \end{bmatrix}$$
% In this case 
% $$\pi(x) = \begin{cases}
%     \frac{1}{4},x=0\\
%     \frac{3}{4}\cdot \frac{1}{2^k}
% \end{cases}$$
$$h_i(x) = \begin{cases}
    0,x\leq a-1\\
    \frac{(a+1)^2}{3},x=a
\end{cases},\,i=1,2.$$

% So we have an accurate solution $\pi r=\frac{3}{4}$
 We use the Python package ``np.linalg.solve" to solve linear equations. For the algorithm in \cite{infanger2022new}, there is one step that uses the inverse function in the Julia language, which does not take advantage of sparsity. For our comparison, we replaced the inverse operation by instead solving systems of linear equations. The comparison results can be found in Table 1; Reg is the algorithms introduced in this paper, whereas Trunc is the algorithm studied in \cite{infanger2022new}. Table 1 provides a comparison of computation times and error bound quality. The results suggest that the two algorithms have comparable run times and error bound quality.
 
    \begin{center}
    \begin{table}[h]
        \centering
        \begin{tabular}{lcccccc}
            \toprule
            \multicolumn{4}{c}{G/M/1 Model} & \multicolumn{3}{c}{Random Walk}\\
            \toprule
            ~ &  $a = 10^3$  & $a = 5 \times 10^3$ & $a = 10^4$ & $a = 10^3$  & $a = 5 \times 10^3$ & $a = 10^4$\\
            \midrule

            \multicolumn{7}{c}{CPU time}\\

            Reg  & 0.24   &  6.77  & 43.45 & 0.418& 6.39 & 45.49  \\
            Trunc  & 0.35 & 7.48   & 44.12 & 0.384& 7.26 & 45.44 \\
                        \hline
            \multicolumn{7}{c}{Upper Bounds}\\
            Reg  & 137.548   & 133.167   & 133.167   & 0.75000 & 0.75000& 0.75000   \\
            Trunc  & 137.348 & 133.167   & 133.167   & 0.74999 & 0.74999& 0.74999  \\ 
        \hline
            \multicolumn{7}{c}{Lower Bounds}\\

            Reg  & 130.147   & 133.167    & 133.167   & 0.74999& 0.74999  & 0.74999   \\
            Trunc & 18.02 & 133.167    & 133.167    &0.74999 &0.74999  & 0.74999   \\
            \hline         
            \bottomrule
        \end{tabular}
        \caption{Comparison for computation time}
    \end{table}
\end{center}

%  \begin{center}
%     \begin{table}[h]
%         \centering
%         \begin{tabular}{lcccccc}
%             \toprule
%             \multicolumn{4}{c}{G/M/1 Model} & \multicolumn{3}{c}{Random Walk}\\
%             \toprule
%             ~ &  $A = 10^3$  & $A = 5 \times 10^3$ & $A = 10^4$ & $A = 10^3$  & $A = 5 \times 10^3$ & $A = 10^4$\\
%             \midrule

%             \multicolumn{7}{c}{CPU time}\\

%             Reg  & 0.24   &  6.77  & 43.45 & 0.418& 6.39 & 45.49  \\
%             Trunc  & 0.35 & 7.48   & 44.12 & 0.384& 7.26 & 45.44 \\
%             \hline
%             \multicolumn{7}{c}{Upper Bounds}\\
%             Reg  & 137.548   & 133.167   & 133.167   & 0.75000 & 0.75000& 0.75000   \\
%             Trunc  & 137.348 & 133.167   & 133.167   & 0.74999 & 0.74999& 0.74999  \\ 
%         \hline
%             \multicolumn{7}{c}{Lower Bounds}\\

%             Reg  & 130.147   & 133.167    & 133.167   & 0.74999& 0.74999  & 0.74999   \\
%             Trunc & 18.02 & 133.167    & 133.167    &0.74999 &0.74999  & 0.74999   \\
%             \hline          
%             \bottomrule
%         \end{tabular}
%         \caption{The Comparison of Two Algorithms for the two settings}
%     \end{table}
% \end{center}  

\bibliographystyle{informs2014} % outcomment this and next line in Case 1
\bibliography{ref}

%% Here starts the e-companion (EC)
%%%%%%%%%%%%%%%%%%%%%%%%%%%%%%%%%%%%%%%%%%%%%%%%%%%%%%%%%%
% \ECSwitch \small

% \newpage

% \section*{Appendix} \label{appendix}
% %\ECDisclaimer
% %%%%%%%%%%%%%%%%%%%%%%%%%%%%%%%%%%%%%%%%%%%%%%%%%%%%%%%%%%

% %%% Main head for the e-companion
% % \ECHead{Proofs of Statements}

% \subsection*{Proof of Theorem 1} \label{appendixmse}

% Acknowledgments here
% \ACKNOWLEDGMENT{The authors gratefully acknowledge...}

\end{document}